\documentclass[12pt,reqno]{amsart}
\usepackage[utf8]{inputenc}

\usepackage{amsfonts, amssymb, amsmath, amsthm, mathtools}
\usepackage{mathrsfs} 
\usepackage{latexsym}
\usepackage{enumerate}
\usepackage{multicol}
\usepackage{times}
\usepackage{verbatim}
\usepackage{tikz-cd}
\usepackage{fullpage}
\usepackage{here}
\usepackage{url}
\usepackage{csquotes}
\usepackage{placeins}
\usepackage{epic}
\usepackage{graphicx}
\usepackage{epstopdf}
\usepackage{pstricks}
\usepackage{pst-plot}
\usepackage{tikz}
\usetikzlibrary{shapes.geometric, positioning, calc}

\input xypic
\xyoption{all}
\xyoption{poly}

\usepackage[colorlinks=true]{hyperref}
\hypersetup{citecolor=blue, linkcolor=blue}

\usepackage[capitalise, nameinlink]{cleveref}

\crefname{equation}{}{}
\crefname{figure}{{\sc Figure}}{{\sc Figure}}
\crefname{subsection}{Subsection}{Subsections}

\usetikzlibrary{matrix,arrows,decorations.pathmorphing}

\usepackage[euler]{textgreek}
\usepackage{tikz,tikz-cd}
\usepackage[colorinlistoftodos, textwidth = 2.3cm]{todonotes}

\newtheorem{theorem}{Theorem}[section]

\newtheorem{lemma}[theorem]{Lemma}
\newtheorem{corollary}[theorem]{Corollary}

\newtheorem*{claim*}{Claim}

\theoremstyle{definition}

\newtheorem{remark}[theorem]{Remark}

\newcommand{\F}{{\mathbb F}}
\newcommand{\Z}{{\mathbb Z}}

\DeclareMathOperator{\ord}{ord}
\usepackage{tabstackengine}
\stackMath

\numberwithin{equation}{section} 
\numberwithin{figure}{section}
\numberwithin{table}{section}

\title{$F$-Diophantine sets over finite fields}

\author{Chi Hoi Yip}
\address{School of Mathematics\\ Georgia Institute of Technology\\ Atlanta, GA 30332\\ United States}
\email{cyip30@gatech.edu}
\author{Semin Yoo}
\address{Discrete Mathematics Group \\ Institute for Basic Science \\ 55 Expo-ro Yuseong-gu, Daejeon 34126 \\ South Korea}
\email{syoo19@ibs.re.kr}
\subjclass[2020]{11D79, 11T06, 11T24}
\keywords{Diophantine tuple, $F$-Diophantine set, finite field}
\begin{document}

\begin{abstract}
Let $k \geq 2$, $q$ be an odd prime power, and $F \in \F_q[x_1, \ldots, x_k]$ be a polynomial. An $F$-Diophantine set over a finite field $\F_q$ is a set $A \subset \F_q^*$ such that $F(a_1, a_2, \ldots, a_k)$ is a square in $\F_q$ whenever $a_1, a_2, \ldots, a_k$ are distinct elements in $A$. In this paper, we provide a strategy to construct a large $F$-Diophantine set, provided that $F$ has a nice property in terms of its monomial expansion. In particular, when $F=x_1x_2\ldots x_k+1$, our construction gives a $k$-Diophantine tuple over $\F_q$ with size $\gg_k \log q$, significantly improving the $\Theta((\log q)^{1/(k-1)})$ lower bound in a recent paper by Hammonds-Kim-Miller-Nigam-Onghai-Saikia-Sharma.
\end{abstract}

\maketitle

\section{Introduction}
A set of $m$ positive integers is a \textit{Diophantine $m$-tuple} if the product of any two distinct elements in the set is one less than a perfect square. There are many interesting results in the study of Diophantine tuples and their variants. Perhaps most notable is the Diophantine quintuple conjecture, namely, there is no Diophantine quintuple, recently confirmed by He, Togb\'e, and Ziegler \cite{HTZ19}. We refer to the Dujella's book \cite{D24} for a comprehensive discussion on the topic and their reference.

The definition of $F$-Diophantine sets were formally introduced by B\'{e}rczes, Dujella, Hajdu, Tengely \cite{BDHT16} for a polynomial $F \in \Z[x,y]$. Given a polynomial $F \in \Z[x,y]$, they say that a subset $A$ of integers is an \emph{$F$-Diophantine set} if $F(x,y)$ is a perfect square for all $x,y \in A$ with $x \neq y$. $F$-Diophantine sets naturally appear in various contexts and are related to many interesting problems in number theory. In particular, an $F$-Diophantine set with $F(x,y)=xy+n$ and $n \neq 0$ corresponds to a Diophantine tuple with property $D(n)$ (see for example \cite{D02}). Similar to the study of classical Diophantine tuples, it is of special interest to construct large $F$-Diophantine sets or give bounds on the maximum size of $F$-Diophantine sets \cite{BDHT16, SE18}. 

In this paper, we study the natural analogue of $F$-Diophantine sets over finite fields. Throughout the paper, let $q$ be an odd prime power, $\F_q$ the finite field with $q$ elements, and $\F_q^*=\F_q \setminus \{0\}$. Let $k \geq 2$ and $F \in \F_q[x_1, \ldots, x_k]$ be a polynomial. We say $A \subset \F_q^*$ is an \emph{$F$-Diophantine set over $\F_q$} if $F(a_1, a_2, \ldots, a_k)$ is a square in $\F_q$ whenever $a_1, a_2, \ldots, a_k$ are distinct elements in $A$. In the same spirit, we are interested in estimating the quantity $M(F;\F_q)$, the maximum size of $F$-Diophantine sets over $\F_q$\footnote{The definition of $M(F;\F_q)$ still makes sense when $q$ is even, however in that case we trivially have $M(F;\F_q)=q-1$ since each element in $\F_q$ is a square.}. Although such terminology appears to be new in general, for many special polynomials $F$, $F$-Diophantine sets over finite fields have been studied extensively in different contexts. The obvious choice $F(x,y)=xy+\lambda$ with $\lambda \in \F_q^*$ corresponds to generalized Diophantine tuples over $\F_q$ \cite{21DK, G01, KYY, KYY24a, Shparlinski23, Y24}. $F$-Diophantine sets over $\F_q$ with $F(x,y)=x-y$ (when $q \equiv 1 \pmod 4$) corresponds to cliques in the Paley graph over $\F_q$. In the aforementioned two cases, when $q$ is a non-square, we have the ``trivial" bounds
$$
(1-o(1))\log_4 q \leq M(F; \F_q) \leq \sqrt{q}+O(1);
$$
see \cite{KYY, KYY24b, Yip1}. However, any bound beyond the above requires highly non-trivial efforts. We refer to \cite{88Cohen, KYY, KYY24a, KYY24b, Yip1, Y23++} for recent multiplicative constant improvement on the lower bounds and upper bounds from polynomial methods, finite geometry, number theory, and graph theory. Moreover, when $k=2$, the authors \cite{KYY24b} have studied lower bounds and upper bounds on $M(F; \F_q)$ for a generic polynomial $F\in \F_q[x,y]$. We focus on the case $k \geq 3$ in this paper.

Next we discuss lower bounds and upper bounds on $M(F;\F_q)$ for a generic polynomial $F \in \F_q[x_1,x_2,\ldots, x_k]$ with degree $d$. From \cite[Section 3.3]{KYY24b}, one can deduce that $M(F;\F_q)=O_{d}(\sqrt{q})$ if $F$ is generic. Note that this upper bound is sometimes sharp. Indeed, if $q$ is a square and $F$ is defined over $\F_{\sqrt{q}}$, then $A=\F_{\sqrt{q}}^*$ is an $F$-Diophantine set over $\F_q$ since all elements in $\F_{\sqrt{q}}$ are squares in $\F_q$. Regarding the lower bound on $M(F; \F_q)$, it is helpful to use a probabilistic heuristic. Assuming that the set of squares in $\F_q$ was a random subset of $\F_q$ with density $1/2$, then we expect that there exists an $F$-Diophantine set over $\F_q$ with size $n$ provided that
$$
\binom{q}{n} 2^{-\binom{n}{k}}\geq 1.
$$
This suggests the heuristic lower bound that
\begin{equation}\label{eq:lb}
M(F;\F_q) \geq \Theta((\log q)^{1/(k-1)}).    
\end{equation}
Here, for two functions $f$ and $g$, $f=\Theta(g)$ means that both $f=O(g)$ and $g=O(f)$ are satisfied.
Indeed, when $F(x_1,x_2,\ldots, x_k)=x_1x_2\cdots x_k+1$, Hammonds, Kim, Miller, Nigam, Onghai, Saikia, and Sharma \cite[Theorem 1.3]{k23} confirmed inequality~\eqref{eq:lb} (they only considered the case where $q$ is an odd prime, but the same proof extends to all odd prime powers $q$). Unsurprisingly, in their terminology, such an $F$-Diophantine set over $\F_q$ is a \emph{$k$-Diophantine tuple over $\F_q$}. 

Before stating our main result, we need to introduce a new definition. We define a partial order on non-constant monic monomials in $\F_q[x_1,x_2, \ldots, x_k]$. Let $f=x_1^{\alpha_1} x_2^{\alpha_2} \cdots x_k^{\alpha_k}$ and $g=x_1^{\beta_1} x_2^{\beta_2} \cdots x_k^{\beta_k}$, where $\alpha_1, \beta_1, \ldots, \alpha_k, \beta_k$ are nonnegative integers. We write $f\succeq g$ if $\alpha_i \geq \beta_i$ for each $1 \leq i \leq k$, and write $f\succ g$ if $f \succeq g$ and $f \neq g$. Let $F \in \F_q[x_1,x_2,\ldots, x_k]$ be a nonzero polynomial. We can write $F$ in its monomial expansion as follows:
$$
F=\sum_{i=1}^m a_i f_i +C,
$$
where each $a_i \in \F_q^*$, $f_i$ is a monomial of degree at least $1$, and $C \in \F_q$. We say $F$ is \emph{of type I} if $C$ is a non-zero square in $\F_q$. We say $F$ is \emph{of type II} if there is $1 \leq i \leq m$, such that $a_i$ is a square in $\F_q^*$, and $f_i \succ f_j$ for all $1 \leq j \leq m$ with $j \neq i$. 

\begin{theorem}\label{thm:main}
Let $q\geq 257$ be an odd prime power and let $F \in \F_q[x_1,x_2,\ldots, x_k]$ be a nonzero polynomial of type I or type II. If $F$ has degree $d$ and the monomial expansion of $F$ consists of $m$ non-constant monomials, then 
$$
M(F;\F_q)\geq \bigg\lfloor\frac{1}{d}(\log_4 q-4\log_4 \log_4 q)^{1/m} \bigg\rfloor.
$$
\end{theorem}

Applying Theorem~\ref{thm:main} to $F(x_1,x_2,\ldots, x_k)=x_1x_2\ldots x_k+1$ (which is of both type I and type II), we get the following corollary immediately. In particular, it significantly improves the lower bound $\Theta((\log q)^{1/(k-1)})$ on the maximum size of $k$-Diophantine tuples over $\F_q$ by Hammonds et. al \cite{k23}.

\begin{corollary}\label{cor:main}
Let $k \geq 2$ and let $q$ be an odd prime power. There is an $k$-Diophantine tuple over $\F_q$ with size at least $\big(\frac{1}{k}-o(1)\big)\log_4 q$, as $q \to \infty$. 
\end{corollary}

Interestingly, our approach provides a substantial improvement on the heuristic lower bound of $M(F;\F_q)$ given in inequality~\eqref{eq:lb}, whenever $k \geq 3$ and $F\in \F_q[x_1,x_2,\ldots, x_k]$ is a sparse polynomial. The constant factors in Theorem~\ref{thm:main} and Corollary~\ref{cor:main} are not optimal. Here we focus on improving the order of the magnitude of the lower bound and we do not attempt to optimize the constant factors. On the other hand, in the case $k=2$, improving the constant factors in front of $\log q$ is of special interest; see \cite{KYY24b} and references therein for more discussions.

\section{Constructions of $F$-Diophantine sets}
Let $q$ be an odd prime power. Let $F \in \F_q[x_1,x_2,\ldots, x_k]$ be a polynomial of type I or type II, with degree $d$. Write $F$ in its monomial expansion as follows:
$$
F=\sum_{i=1}^m a_i f_i +C,
$$
where $a_i \in \F_q^*$ and $f_i$ is a monomial of degree at least $1$ for each $1\leq i \leq m$, and $C \in \F_q$. 

Let $n$ be a positive integer to be determined such that $2 \leq n \leq q^{1/4}$.
Consider the following collection of polynomials in $\F_q[x]$:
$$
V:=V(n)=\{F(x^{\theta_1}, x^{\theta_2}, \ldots, x^{\theta_k}): 1 \leq \theta_1, \theta_2, \ldots, \theta_k \leq n\}.
$$
Observe that
\begin{equation}\label{eq:V}
V \subset \bigg\{\sum_{i=1}^m a_i x^{\alpha_i} +C: 1\leq \alpha_1, \alpha_2, \ldots, \alpha_m \leq dn\bigg\}.    
\end{equation}
Also, if $F$ is of type I, then the constant term of each polynomial $g \in V$ is a non-zero square in $\F_q$; if $F$ is of type II, then the leading coefficient of each polynomial $g \in V$ is a non-zero square in $\F_q$. In both cases, it readily follows that the product of polynomials in any subset of $V$ is not of the form $ch^2$, where $c$ is a non-square in $\F_q$, and $h$ is a polynomial in $\F_q[x]$. 

Let $Y$ denote the collection of $y \in \F_q^*$ with order at least $n$, such that the set $\{g(y): g \in V\}$ is contained in the set of squares in $\F_q$. Let $N=|Y|$ and let $\chi$ be the quadratic character in $\F_q$. We claim that 
\begin{equation}\label{eq:N}
N\geq 2^{-|V|} \sum_{\substack{y\in\F_q^* \\ \ord y \geq n }} \prod_{g \in V}\bigg(1+\chi \big(g(y)\big)\bigg).    
\end{equation}
Indeed, if $y \notin Y$, then $g(y)$ is a non-square in $\F_q$ for some $g \in V$, and thus such $y$ does not contribute to the right-hand side of inequality~\eqref{eq:N}. On the other hand, if $y \in Y$, then $\chi(g(y)) \in \{0,1\}$ for each $g \in V$, and thus it contributes at most $1$ to the right-hand side of inequality~\eqref{eq:N}.

Expanding the product on the right-hand side of inequality~\eqref{eq:N} yields
\begin{align*}
N
&\geq 2^{-|V|} \sum_{\substack{y\in\F_q^* \\ \ord y \geq n }} \sum_{W \subset V} \prod_{g \in W} \chi \big(g(y)\big)\\
&= 2^{-|V|} \sum_{W \subset V} \sum_{\substack{y\in\F_q^* \\ \ord y \geq n }}   \chi \bigg(\bigg(\prod_{g \in W}g\bigg)(y)\bigg)\\
&\geq  -|Z|+ 2^{-|V|} \sum_{W \subset V} \sum_{y \in \F_q}   \chi \bigg(\bigg(\prod_{g \in W}g\bigg)(y)\bigg),
\end{align*}
where $Z=\{0\} \cup \{y \in \F_q^*: \ord y<n\}$. Since $\F_q^*$ is a cyclic group, it is clear that $|Z|\leq n^2$. 

We need to use Weil's bound for complete character sums (see for example \cite[Theorem 5.41]{LN97}), which we recall below. 

\begin{lemma}(Weil's bound)\label{Weil}
Let $\chi$ be a multiplicative character of $\F_q$ of order $k>1$, and let $g \in \F_q[x]$ be a monic polynomial of positive degree that is not an $k$-th power of a polynomial. 
Let $s$ be the number of distinct roots of $g$ in its
splitting field over $\F_q$. Then for any $a \in \F_q$,
$$\bigg |\sum_{x\in\mathbb{F}_q}\chi\big(ag(x)\big) \bigg|\le(s-1)\sqrt q.$$
\end{lemma}

We have mentioned that for each subset $W$ of $V$, the product $\prod_{g \in W}g$ is not of the form $ch^2$, where $c$ is a non-square in $\F_q$, and $h$ is a polynomial in $\F_q[x]$. Therefore, separating the contribution from $W=\emptyset$ and $W \neq \emptyset$, and applying Weil's bound, we further deduce that
\begin{align}
N&\geq -n^2+ 2^{-|V|} \sum_{W \subset V} \sum_{y \in \F_q}   \chi \bigg(\bigg(\prod_{g \in W}g\bigg)(y)\bigg) \notag\\
&\geq -n^2+\frac{q}{2^{|V|}}- 2^{-|V|}\sum_{\substack{W \subset V\\W \neq \emptyset}} \bigg(-1+\sum_{g \in W} \deg(g)\bigg)\sqrt{q} \notag\\
&=\frac{q}{2^{|V|}}-n^2+\frac{\sqrt{q} (2^{|V|}-1)}{2^{|V|}}-2^{-|V|}\sqrt{q}\sum_{W \subset V} \sum_{g \in W} \deg(g). \label{eq:N2}
\end{align}
Given inclusion~\eqref{eq:V}, $\deg(g)\leq dn$ for each $g\in V$. Thus, a simple double-counting argument shows that
\begin{equation}\label{eq:counting}
\sum_{W \subset V} \sum_{g \in W} \deg(g)
=2^{|V|-1}\sum_{g \in V} \deg(g) \leq 2^{|V|-1} |V| dn.   
\end{equation}

We conclude from the assumption $n \leq q^{1/4}$, inequality~\eqref{eq:N2} and inequality~\eqref{eq:counting} that
$$
N \geq \frac{q}{2^{|V|}}-n^2+ \frac{\sqrt{q} (2^{|V|}-1)}{2^{|V|}}-\frac{|V| dn\sqrt{q}}{2} \geq \frac{q}{2^{|V|}}-|V| dn\sqrt{q}.
$$
Note that $|V|\leq (dn)^m$ from  inclusion~\eqref{eq:V}, thus
\begin{equation}\label{eq:N3}
N \geq \frac{q}{2^{(dn)^m}}-(dn)^{m+1}\sqrt{q}.    
\end{equation}
Since $q \geq 257$, we have $\log_4 q>4\log_4 \log_4 q$. Set $$n=\bigg\lfloor\frac{1}{d}(\log_4 q-4\log_4 \log_4 q)^{1/m} \bigg\rfloor.$$ Then we have $(dn)^m\leq \log_4 q-4\log_4 \log_4 q$ and thus 
$$
4^{(dn)^m} (dn)^{2m+2}\leq 4^{(dn)^m} (dn)^{4m}\leq\frac{q}{(\log_4 q)^4} \cdot (dn)^{4m}<q.
$$
It follows from inequality~\eqref{eq:N3} that
$$
N\geq \frac{q}{2^{(dn)^m}}-(dn)^{m+1}\sqrt{q}>0.
$$
Note that $N>0$ implies that $N \geq 1$, that is, there exists $y_0 \in \F_q^*$ with order at least $n$, such that the set $\{g(y_0): g \in V\}$ is contained in the set of squares in $\F_q$. Let 
$$
A=\{y_0^1, y_0^2, \ldots, y_0^n\};
$$
then $A$ is an $F$-Diophantine set over $\F_q$ with $|A|=n$. This proves Theorem~\ref{thm:main}, as required.

\medskip

Next, we give several remarks on our constructions. 

\begin{remark}
Our constructions above in fact produce a \emph{strong $F$-Diophantine set} $A$ over $\F_q$ in the sense that $F(a_1, a_2, \ldots, a_k)$ is a square in $\F_q$ whenever $a_1, a_2, \ldots, a_k$ are elements in $A$ (not necessarily distinct), in the spirit of strong Diophantine tuples \cite{DP08, KYY}. In many cases, one can modify the definition of $V$ in the above construction to obtain a slightly larger $F$-Diophantine set over $\F_q$. 
\end{remark}

\begin{remark}
In general, to construct a large $F$-Diophantine set over $\F_q$, we need to impose some assumptions on the polynomial $F$. Obviously, we have to assume that $F$ is not of the form $cG^2$, where $c$ is a non-square in $\F_q$ and $G \in \F_q[x_1,x_2, \ldots, x_k]$. 

The assumption that $F$ is of type I or type II made in the statement of Theorem~\ref{thm:main} can be weakened. Indeed, as long as one can come up with a similar definition of $V$, and show that the product of polynomials in any subset of $V$ is not of the form $ch^2$ (where $c$ is a non-square in $\F_q$, and $h$ is a polynomial in $\F_q[x]$), then one can modify the above proof to produce a large $F$-Diophantine set over $\F_q$.

As an illustration, consider a degree $d$ homogeneous polynomial $$F(x_1,x_2,\ldots, x_k)=\sum_{i=1}^k c_ix_i^d \in \F_q[x_1,x_2,\ldots, x_k]$$ with $k \geq 2$, where $c_i$ is a non-zero square in $\F_q$ for each $1 \leq i \leq k$. Note that such $F$ is neither of type I nor type II. If we instead define
$$
V:=V(n)=\{F(x^{\theta_1}, x^{\theta_2}, \ldots, x^{\theta_k}): \theta_1, \theta_2, \ldots, \theta_k \text{ are distinct elements in } \{1,2,\ldots, n\}\},
$$
then the leading coefficient of each polynomial $g \in V$ is a non-zero square in $\F_q$. It follows that the product of polynomials in any subset of $V$ is not of the form $ch^2$, where $c$ is a non-square in $\F_q$, and $h$ is a polynomial in $\F_q[x]$. Thus, a similar argument as above shows that there is an $F$-Diophantine set $A$ over $\F_q$ with $|A|=n$, where
$$
n \geq \bigg(\frac{1}{d}-o(1)\bigg) (\log_{4} q)^{1/k}.
$$
Note however in this case, if $\sum_{i=1}^k c_i$ is a non-square in $\F_q$, and $d$ is even, then there is no strong $F$-Diophantine set over $\F_q$.
\end{remark}

\begin{remark}
Recently, there have been a few papers devoted to the search for Diophantine tuples with additional properties. For example, looking for a Diophantine tuple with property $D(n)$ for multiple different $n$ \cite{CGH23}, or a rational Diophantine tuple with square elements \cite{DKP21}. In a similar flavor, it would be interesting to search for a large $F$-Diophantine set over a finite field with additional properties, and usually it is not hard to modify the above proof to achieve this purpose. For example, if we want to look for a large $F$-Diophantine set over $\F_q$ with size $n$ consisting of square elements, we can simply change the definition of $V$ to 
$$
V:=V(n)=\{F(x^{2\theta_1}, x^{2\theta_2}, \ldots, x^{2\theta_k}): 1 \leq \theta_1, \theta_2, \ldots, \theta_k \leq n\}
$$
and modify the above proof accordingly. 
\end{remark}

\section*{Acknowledgements}
The second author was supported by the Institute for Basic Science (IBS-R029-C1). The authors thank Seoyoung Kim for helpful discussions. The authors are also grateful to anonymous referees for their valuable comments and corrections. 

\bibliographystyle{abbrv}
\bibliography{references}

\end{document}